\title{ ~~\\ On the distribution of the order and index of
$g({\rm mod~}p)$ over residue classes III}
\author{Pieter Moree}
\documentclass[12pt]{article}
\usepackage{amssymb, latexsym, amsfonts}
\def\@ptsize{2}
\newtheorem{Thm}{Theorem}

\newtheorem{Lem}{Lemma}
\newtheorem{Def}{Definition}
\newtheorem{Cor}{Corollary}

\newtheorem{Prop}{Proposition}
\newcommand{\qed}{\hfill $\Box$}

\begin{document}
\date{}
\maketitle
{\def\thefootnote{}
\footnote{\noindent 
Max-Planck-Institut f\"ur Mathematik, Vivatsgasse 7,
D-53111 Bonn, Deutschland, E-mail:
moree@mpim-bonn.mpg.de\\}
{\def\thefootnote{}
\footnote{{\it Mathematics Subject Classification (2000)}.
11N37, 11N69, 11R45}}

\begin{abstract}
\noindent For a fixed rational number $g\not\in \{-1,0,1\}$ and
integers $a$ and $d$ we consider the
sets  $N_g(a,d)$, respectively $R_g(a,d)$, of primes $p$ for which the order, 
respectively the index of $g({\rm mod~}p)$ is
congruent to $a({\rm mod~}d)$.
Under the
Generalized Riemann Hypothesis (GRH), it is known that these sets have a natural
density $\delta_g(a,d)$, respectively $\rho_g(a,d)$. 
It is shown that these densities can be expressed as linear combinations of 
certain constants introduced by Pappalardi.
Furthermore it is proved that
 $\delta_g(a,d)$ and $\rho_g(a,d)$ equal their $g$-averages for almost all $g$.
 \end{abstract}

\section{Introduction}
\label{introdu}
Let $g\not\in \{-1,0,1\}$ be a rational number (this 
assumption on $g$ will be maintained throughout this paper) and let $a$ and $d$ be
positive integers. The central object of study of this paper and its predecessors 
\cite{Moreealleen, Moreealleen2} are the sets $N_g(a,d)$ and, 
to a lesser extent, due to its much greater 
simplicity of description, $R_g(a,d)$. Nevertheless, the main results for both sets 
(Theorems \ref{stellingeen}, \ref{favoriet}, \ref{indeks}) are quite
analogous in character. 
However, Theorem \ref{peelof} does not have an analog for $\rho$, cf. Proposition \ref{nietanaloog}.
The set $N_g(a,d)$ consists of the primes $p$ for which
$\nu_p(g)=0$ and ord$_g(p)\equiv a({\rm mod~}d)$ (here and
in the sequel the letter $p$ will be used to indicate prime 
numbers). The set $R_g(a,d)$ is similarly defined, but with ord$_g(p)$ replaced
by $r_g(p)$, the residual index.\\
\indent Pappalardi \cite{Pappalardi} established, under GRH, the
existence of the natural density $\rho_g(a,d)$ of $R_g(a,d)$. In particular, he showed
that
\begin{equation}
\label{pappastelling}
\rho_g(a,d)=\sum_{t\equiv a({\rm mod~}d)}\sum_{n=1}^{\infty}
{\mu(n)\over [K_{nt,nt}:\mathbb Q]},
\end{equation}
where for $r|s$ we define $K_{s,r}$ to be the number field $\mathbb Q(\zeta_s,g^{1/r})$.
In \cite{Moreealleen2} it was shown that, under GRH, the set $N_g(a,d)$ has a 
natural density $\delta_g(a,d)$. In particular, 
\begin{equation}
\label{centraal}
\delta_g(a,d)=
\sum_{t=1\atop (1+ta,d)=1}^{\infty}
\sum_{n=1\atop (n,d)|a}^{\infty}{\mu(n)c_g(1+ta,dt,nt)\over 
[K_{[d,n]t,nt}:\mathbb Q]},
\end{equation}
where, for $(b,f)=1$,
$$c_g(b,f,v)=\cases{1 &if $\sigma_b|_{\mathbb Q(\zeta_{f})\cap K_{v,v}}={\rm id}$;\cr
0 & otherwise,}$$
where $\sigma_b$ is the automorphism of $\mathbb Q(\zeta_f)$ that sends
$\zeta_f$ to $\zeta_f^b$
(if
$a$ and $b$ are integers, then by $(a,b)$ and $[a,b]$ we denote
the greatest common divisor, respectively lowest common multiple of
$a$ and $b$).
The rational number $c_g(1+ta,dt,nt)/[K_{[d,n]t,nt}:\mathbb Q]$ is
the density of primes $p\equiv 1+ta({\rm mod~}dt)$ that split completely in $K_{nt,nt}$. It turns
out that the coefficient $c_g(.,.,.)$ has a strong tendency to equal $1$ and this motivates the
following definition:
\begin{Def}
Put
\begin{equation}
\label{centraal1}
\delta_g^{(0)}(a,d)=
\sum_{t=1\atop (1+ta,d)=1}^{\infty}
\sum_{n=1\atop (n,d)|a}^{\infty}{\mu(n)\over 
[K_{[d,n]t,nt}:\mathbb Q]},
\end{equation}
\end{Def}
For example, it can be shown that (unconditionally) $\delta_g(0,d)=\delta_g^{(0)}(0,d)$. Other
examples are provided by Proposition \ref{other} and Lemma \ref{bijnagelijktwee}.\\
\indent In
case $d=2$ unconditional results can be obtained, cf. \cite{Odoni, Wiertelak1, Wiertelak2}.
For an extensive analysis of the case $d=3$, see \cite{Moreealleen}. For the
case $d=4$ see also \cite{Moreealleen} (for less general results obtained by a
different method for this modulus see \cite{CM}). The case where $d$ is prime
is investigated in \cite{Moreealleen2}. In this paper the general case is investigated. Not so
surprisingly the present results will be somewhat less explicit.\\
\indent  By the following result $d$ can
be taken close to squarefree. 
\begin{Thm} {\rm (GRH)}.
\label{peelof}$~$\\
{\rm 1)} If $q$ is an odd prime dividing $d_1$, then $\delta_g(a,qd_1)=\delta_g(a,d_1)/q$.\\
{\rm 2)} If $8|d$, then $\delta_g(a,2d_1)=\delta_g(a,d_1)/2$.\\
The same conclusion holds with $\delta_g$ replaced by $\delta_g^{(0)}$ or $\delta$. In the latter
case, the result holds true unconditionally.
\end{Thm}
By $k(d)$ we denote the 
squarefree part of $d$, that is $k(d)=\prod_{p|d}p$.
Also we define 
$$k_1(d)=\cases{k(d) & if $d$ is odd;\cr
4k(d) & otherwise,}{\rm ~and~}k_2(d)=\cases{k(d) & if $d$ is odd;\cr
(4,d/2)k(d) & otherwise.}$$
Note that, assuming GRH, we have $\rho_g(a,d)=\rho_g(a,d_1)d_1/d$, where $d_1=k_2(d)$.\\
\indent It turns out that the average density of elements in a finite field of prime 
cardinality having order, respectively
index, congruent to $a({\rm mod~}d)$ exists. Denote these densities by $\delta(a,d)$, respectively $\rho(a,d)$.
These densities can be studied by quite elementary methods and are easier to study than $\delta_g(a,d)$ and
$\rho_g(a,d)$. In this paper we will be particularly interested in the connection between $\delta_g(a,d)$ 
and $\delta(a,d)$ and the same for $\rho_g(a,d)$ and $\rho(a,d)$. The quantities $\delta(a,d)$ and
$\rho(a,d)$ are studied in \cite{Moreeaverage}. In this paper we will show that, under GRH, $*(a,d)$ can
also be regarded as
the $g$-average of $*_g(a,d)$ (where for ease of notation we define $\delta^{(0)}(a,d)=\delta(a,d)$): 
\begin{Prop} {\rm (GRH)}.
\label{gaverage}
With $*=\delta,\delta^{(0)}$ or $\rho$ we have
$${1\over 2x}\sum_{|g|\le x}*_g(a,d)=*(a,d)+O({1\over \sqrt{x}}).$$
\end{Prop}
(All results in this paper with the condition $|g|\le x$ remain valid if this condition is replaced
with $1<g\le x$ or $-x\le g<1$ and $2x$ by $x$.) Proposition \ref{gaverage} is a simple consequence
of the following result in which $G$ is the set of rational numbers $g$ that cannot be written as
$-g_0^h$ or $g_0^h$ with $h>1$ an integer and $g_0$ a rational number (in other words the rational
numbers $g$ that satisfy $h=1$ in the notation of Lemma \ref{degree}).
\begin{Thm} {\rm (GRH)}. 
\label{stellingeen}
Suppose that $g\in G$. Set $D_1=|D(g)/(D(g),d)|$ and
$D_2=[2,D(g)]$, where $D(g)$ denotes the discriminant of the number field $\mathbb Q(\sqrt{g})$.
Then
$$\Big|\delta_g^o(a,d)-\delta(a,d)\Big|< {2^{\omega(D_1)+2}\over \varphi(D_1)D_1},~~~\Big|\delta_g(a,d)-\delta(a,d)\Big|<{3\cdot 2^{\omega(D_1)+2}\over \varphi(D_1)D_1},$$
and $$\Big|\rho_g(a,d)-\rho(a,d)\Big|< {2^{\omega(D_2)+2}\over \varphi(D_2)D_2},$$
where $\omega(n)$ denotes the number of distinct prime divisors of $n$.
\end{Thm}
\begin{Cor}
\label{gggevolg}
If $|D(g)|$ tends to infinity with $g\in G$, then $*_g(a,d)$ tends to $*(a,d)$, 
where $*=\delta,\delta^{(0)}$ or $\rho$.
\end{Cor}
Corollary \ref{gggevolg} for $*=\delta_g$ has been already announced (but not proved) as Theorem 2 of \cite{Moreeaverage}.
Theorem \ref{stellingeen} shows that `generically' $*_g(a,d)$ is quite close to $*(a,d)$.
For a numerical demonstration of this see Tables 1 and 2 of \cite{Moreeaverage}. Much more is true, however:
\begin{Thm} {\rm (GRH)}.
\label{favoriet}
Let $d$ be fixed.
There are at most $O_d(x\log^{-1/\varphi(k_1(d))}x)$ integers $g$ with $|g|\le x$ such that $*_g(a,d)\ne *(a,d)$ for
some integer $a$, 
where $*=\delta,\delta^{(0)}$ or $*=\rho$. In particular,
$$(*_g(0,d),\dots,*_g(d-1,d))=(*(0,d),\dots,*(d-1,d))$$
for almost all integers $g$.
\end{Thm}
This result allows one to introduce various notions of `genericity' for $g$.
Actually a more precise version of Theorem \ref{favoriet} is provided 
by Theorem \ref{dgisd} in case $*=\delta$ or $*=\delta^{(0)}$ and
by Propositions \ref{rhogisrho} and \ref{scherper}
in case $*=\rho$. Note that $k_1(d)|4d$.\\
\indent The final main result is concerned with explicitly evaluating $*_g(a,d)$.
For $\chi$ a Dirichlet character, let 
$$A_{\chi}=\prod_{p\atop \chi(p)\ne 0}\left(1+{[\chi(p)-1]p\over [p^2-\chi(p)](p-1)}\right).$$
The constants $A_{\chi}$ turn out to be the basic constants in this
context. They were introduced by Pappalardi \cite{Pappalardi}.
In many cases $A_{\chi}\in \mathbb C\backslash \mathbb R$, see 
\cite[Table 3]{Moreeaverage}. It can be shown that
$$A_{\chi}\prod_{p|d}(1-{1\over p(p-1)})=A{L(2,\chi)L(3,\chi)\over L(6,\chi^2)}\prod_{r=1}^{\infty}\prod_{k=3r+1}^{\infty}
L(k,\chi^r)^{\lambda(k,r)},$$
where $A$ denotes the Artin constant and the numbers $\lambda(k,r)$ are non-zero integers that
can be related to Fibonacci numbers \cite{PFibo} and $L(s,\chi)$ denotes the Dirichlet L-series 
associated to the character $\chi$. The latter expansion 
of $A_{\chi}$ can be used to approximate $A_{\chi}$ with high numerical accuracy 
(see \cite[Section 6]{Moreeaverage}).\\
\indent
(By $G_d$ we denote the group
of Dirichlet characters mod $d$ and by $o_{\chi}$ the order of the Dirichlet character in $G_{\chi}$.)
\begin{Thm} {\rm (GRH)}.
\label{indeks}
Let $a$ and $d$ be arbitrary natural numbers. Then there exists an integer $d_1$ such that
$$*(a,d)=\sum_{\chi\in G_{d_1}}c_{\chi}A_{\chi}{\rm ~with~}c_{\chi}\in \mathbb Q(\zeta_{o_{\chi}}).$$
Furthermore, $c_{\chi}$ can be explicitly computed.\\
{\rm 1)} If $*=\delta_g^{(0)}$, then one can take $d_1=k_2(d)$.\\
{\rm 2)} If $*=\delta_g$, then one can take $d_1=k_1(d)$.\\
{\rm 3)} If $*=\rho_g$, then one can take $d_1=d/(a,d)$.
\end{Thm}
Unconditional results of this nature for $\delta(a,d)$ and $\rho(a,d)$ are
proved in \cite{Moreeaverage} (with $d_1=k(d)$, respectively $d_1=d/(a,d)$).\\
\indent The proofs of Theorems \ref{favoriet} and \ref{indeks} require more explicit knowledge of the Galois coefficients $c_g$ appearing in (\ref{centraal}). Theorem 
\ref{maindegree} (which requires a good deal of work in order to be proved) allows one to compute
$\mathbb Q(\zeta_f)\cap K_{v,v}$. Using this the coefficient $c_g$ can be explicitly related to
a Kronecker symbol. In all this, unfortunately, the case $g<0$ turns out to be considerably more
complicated than the case $g>0$.

\section{Preliminaries on algebraic number theory}
Let $K:\mathbb Q$ be an abelian number field. By the Kronecker-Weber theorem there exists an
integer $f$ such that $K\subseteq \mathbb Q(\zeta_f)$. The smallest such integer is called the
conductor of $K$. Note that $K\subseteq \mathbb Q(\zeta_n)$ iff $n$ is divisible by the conductor.
Note also that the conductor of a cyclotomic field is never congruent to $2({\rm mod~}4)$.
The following lemma allows one to determine all quadratic subfields of a given cyclotomic field.
\begin{Lem}
\label{conductor}
The conductor of a quadratic number field is equal to the absolute value of its discriminant.
\end{Lem} 
\indent Consider the cyclotomic extension $\mathbb Q(\zeta_f):\mathbb Q$. There are $\varphi(f)$ distinct
automorphisms each determined uniquely by $\sigma_a(\zeta_f)=\zeta_f^a$, with $1\le a\le f$ and
$(a,f)=1$. We need to know when the restriction of such an automorphism to a given quadratic
subfield of $\mathbb Q(\zeta_f)$ is the identity. In this direction we have (for a description of
the Kronecker symbol see e.g. \S 2.1 of \cite{Moreealleen2}):
\begin{Lem}
\label{conductor2}
Let $\mathbb Q(\sqrt{g_1})\subseteq \mathbb Q(\zeta_f)$ be a quadratic field of discriminant $D(g_1)$ and
$a$ be an integer with $(a,f)=1$.
We have $\sigma_a|_{\mathbb Q(\sqrt{g_1})}={\rm id}$ iff $({D(g_1)\over a})=1$, with
$({\cdot \over \cdot})$ the Kronecker symbol.
\end{Lem}

\noindent From the theory of profinite groups we recall the
notion of a supernatural (or Steinitz) number. A supernatural number is
a formal product $\prod_p p^{e_p}$, where each
$e_p\in \mathbb N\cup \{\infty\}$. The set of supernatural numbers
forms a commutative monoid with respect to the obvious product.
If $a$ is a supernatural
number, then by $\nu_p(a)$ we denote the exponent of $p$ 
occurring in $a$. We say that $a|b$, that is
that $a$ divides $b$, if $\nu_p(a)\le \nu_p(b)$ for all primes $p$.
If $a|b$, then $b=ca$, with $c$ a supernatural number. Given supernatural
numbers, we define their greatest common divisor by the formula
$(a,b)=\prod_p p^{{\rm inf}(\nu_p(a),\nu_p(b))}$. To a natural number
$d$ we associate the supernatural number $d_{\infty}$, where
$d_{\infty}=\prod_{p|d}p^{\infty}$. The `$d$-part' of a natural
number $n$ is then given by $(n,d_{\infty})$.

\subsection{Preliminaries on field degrees and intersections}
In order to explicitly evaluate certain densities in this paper, the 
following result will play a crucial role. 
Let $g_1\ne 0$ be a rational number. By $D(g_1)$ we denote the discriminant
of the field $\mathbb Q(\sqrt{g})$. An integer $D$ is said to be a fundamental
discriminant if $D=D(g_1)$ for some $g_1\in \mathbb Q$. If $D$ is a fundamental
discriminant, then $\mathbb Q(\sqrt{D})$ has conductor $|D|$ by Lemma \ref{conductor}.
The notation $D(g_1)$ along with the notation $g_0$,  
$h$, $m$ and $n_r$ introduced in the next lemma will reappear again and again in the sequel. 
\begin{Lem} {\rm \cite{Moreealleen}}.
\label{degree}
Write $g=\pm g_0^h$, where $g_0$ is
positive and not an
exact power of a rational. 
Let $D(g_0)$ denote the discriminant of the field $\mathbb Q(\sqrt{g_0})$.
Put $m=D(g_0)/2$ if $\nu_2(h)=0$ and $D(g_0)\equiv 4({\rm mod~}8)$
or $\nu_2(h)=1$ and $D(g_0)\equiv 0({\rm mod~}8)$, and
$m=[2^{\nu_2(h)+2},D(g_0)]$ otherwise. 
Put $$n_r=\cases{m &if $g<0$ and $r$ is odd;\cr
[2^{\nu_2(hr)+1},D(g_0)] &otherwise.}$$
We have
$$[K_{kr,k}:\mathbb Q]=[\mathbb Q(\zeta_{kr},g^{1/k}):\mathbb Q]={\varphi(kr)k\over
\epsilon(kr,k)(k,h)},$$
where, for $g>0$ or $g<0$ and $r$ even we have
$$\epsilon(kr,k)=\cases{2 &if $n_r|kr$;\cr
1 &if $n_r\nmid kr$,}$$
and for $g<0$ and $r$ odd we have
$$\epsilon(kr,k)=\cases{2 &if $n_r|kr$;\cr
 {1\over 2} &if $2|k$ and $2^{\nu_2(h)+1}\nmid k$;\cr
1 &otherwise.}$$
\end{Lem}
{\tt Remark}. Note that if $h$ is odd, then $n_r=[2^{\nu_2(r)+1},D(g)]$. Note that
$n_r=n_{\nu_2(r)}$.\\

\noindent An easy consequence of the latter result is the following lemma.
\begin{Lem}
\label{graadoprekking} $~$\\
{\rm 1)} Suppose that $q$ is an odd prime dividing $d_1$ and $q^2\nmid n$.
Then $[K_{[qd_1,n]t,nt}:\mathbb Q]=q[K_{[d_1,n]t,nt}:\mathbb Q]$.\\
{\rm 2)} Suppose that $\nu_2(d_1)\ge {\rm max}\{2,\nu_2(D(g_0))\}$ and $4\nmid n$. Then 
$[K_{[2^{\alpha}d_1,n]t,nt}:\mathbb Q]=2^{\alpha}[K_{[d_1,n]t,nt}:\mathbb Q]$.
\end{Lem}

\noindent Let $g_1\ne 0$ be a rational number and $n|2m$. Note that the field $\mathbb Q(\zeta_m,\zeta_{n}\sqrt{g_1})$ is
either equal to $\mathbb Q(\zeta_m)$ or a quadratic extension thereof. The following
result records precisely when the extension is not quadratic. In the case $m=n$ this result
is Lemma 3 of \cite{GP} (but earlier references undoubtedly exist). 
\begin{Lem}
\label{degree1}
Let $g_1\ne 0$ be a rational number and suppose that $n|2m$.
Then we have that
$[\mathbb Q(\zeta_{m},\zeta_{n}\sqrt{g_1}):\mathbb Q(\zeta_m)]=\epsilon_1(m,n)$, where
$$\epsilon_1(m,n)=\cases{1 &if $D(g_1)|m$ and $\nu_2(n)\le \nu_2(m)$;\cr
1 &if $D(g_1)|m$ and $2||n$;\cr
1 &if $D(-g_1)|m$ and $4||n$; \cr
1 &if $D(2g_1)|m$ and $8||n$;\cr
2 &otherwise.}
$$
\end{Lem}
{\it Proof}. If $p|D(g_1)$ and $p$ is odd, then $p$ ramifies in $\mathbb 
Q(\zeta_{m},\zeta_{n}\sqrt{g_1})$. This shows that if in addition $p\nmid n$,
then $[\mathbb Q(\zeta_{m},\zeta_{n}\sqrt{g_1}):\mathbb Q(\zeta_m)]=2$ (since
then $p$ does not ramify in $\mathbb Q(\zeta_m)$). The rest of the
argument requires some case distinctions and is left to the interested reader. \qed\\

\noindent {\tt Example}. The result predicts that $\mathbb Q(i,\zeta_8\sqrt{2})=\mathbb Q(i)$.
Indeed, note that $\zeta_8\sqrt{2}=1+i$.

\begin{Lem}
\label{doorisdoor}
Let $L:\mathbb Q$ be abelian with $L\subseteq M$.
If
\begin{equation}
\label{doorsnedegelijk} 
[\mathbb Q(\zeta_f)\cap L:\mathbb Q]=[\mathbb Q(\zeta_f)\cap M:\mathbb Q].
\end{equation}
for every integer $f$, then $L$ is the maximal abelian subfield of $M$.
\end{Lem}
{\it Proof}. Suppose there is an abelian field $L_1$ with $L\varsubsetneq L_1\subset M$.
Let $f(L_1)$ be the conductor of $L_1$. Then
$$[\mathbb Q(\zeta_{f(L_1)})\cap M:\mathbb Q]\ge [L_1:\mathbb Q]>
[L:\mathbb Q]=[\mathbb Q(\zeta_{f(L_1)})\cap L:\mathbb Q].$$
This contradicts equation (\ref{doorsnedegelijk}) with $f=f(L_1)$. \qed\\

\noindent With some additional effort one can explicitly describe $\mathbb Q(\zeta_f)\cap K_{n,n}$. Clearly
this intersection field is abelian and contains $\mathbb Q(\zeta_{(f,n)})$. Let us first compute the absolute
degree of $K_{[f,n],n}$.
We have
\begin{equation}
\label{eerstedoorsnede}
[K_{[f,n],n}:\mathbb Q]={\varphi(f)[K_{n,n}:\mathbb Q]\over [\mathbb Q(\zeta_f)\cap K_{n,n}:\mathbb Q]}.
\end{equation}
On noting that $\varphi((f,n))\varphi([f,n])=\varphi(f)\varphi(n)$, it follows from Lemma \ref{degree} 
and (\ref{eerstedoorsnede}) that
\begin{equation}
\label{tweededoorsnede}
[\mathbb Q(\zeta_f)\cap K_{n,n}:\mathbb Q(\zeta_{(f,n)})]={[\mathbb Q(\zeta_f)\cap K_{n,n}:\mathbb Q]\over \varphi((f,n))}=
{\epsilon([f,n],n)\over \epsilon(n,n)}.
\end{equation}
It is not difficult to infer from Lemma \ref{degree} that the latter quotient is either $1$ or $2$ (so
the apparent possibility $4$ does never arise). We conclude that $\mathbb Q(\zeta_f)\cap K_{n,n}$ is equal
to $\mathbb Q(\zeta_{(f,n)})$ or a quadratic extension thereof.

\begin{Def}
Let $D$ be a fundamental discriminant. Put $b=D/(f,D)$ and
let $f_{\rm odd}$ be the largest odd divisor of $f$.
Put 
$$\gamma(D)=\cases{(-1)^{(f_{\rm odd},D)-1\over 2}(f_{\rm odd},D) & if $\nu_2(f)<\nu_2(D)$;\cr
(-1)^{b-1\over 2}(f,D) & otherwise.}$$
Put $$\gamma_0(D)=\cases{\gamma(D) &if $D\nmid n$ and $D|[f,n]$;\cr 
1 &otherwise.}$$
Put 
$$\gamma_1(D)=\cases{\gamma(D) &if $\nu_2(f)>\nu_2(n)$ and
$D|[f,n]$;\cr
1& otherwise.}$$
\end{Def}
\begin{Lem}
\label{atfundem}
The numbers $\gamma(D)$, $\gamma_0(D)$, $\gamma_1(D)$, ${D\over \gamma_(D)}$, ${D\over \gamma_0(D)}$ and
${D\over \gamma_1(D)}$
are all fundamental discriminants.
\end{Lem}
{\it Proof}. This requires a few case distinctions and is left to the reader. \qed\\

\noindent {\tt Remark}. In this notation the dependence on $f$ and $n$ is suppressed. 
It would be more accurate to write $\gamma_0(D;f,n)$ etcetera. This notation is used in Lemma \ref{DDT}.

\begin{Lem}
\label{oudedoorsnede}
Let $g_1\ne 0$ be a rational number. Then
$$[\mathbb Q(\zeta_f)\cap \mathbb Q(\zeta_n,\sqrt{g_1}):\mathbb Q(\zeta_{(f,n)})]=\cases{2 &if 
$D(g_1)\nmid n$ and $D(g_1)|[f,n]$;\cr
1 & otherwise.}
$$
Furthermore, $$\mathbb Q(\zeta_f)\cap \mathbb Q(\zeta_n,\sqrt{g_1})=\mathbb Q(\zeta_{(f,n)},\sqrt{\gamma_0(D(g_1))}),$$
with $\gamma_0(D(g_1))\ne 1$ iff $D(g_1)\nmid n$ and $D(g_1)|[f,n]$.
\end{Lem}
{\it Proof}. Using Lemma \ref{conductor} we find that 
$$[\mathbb Q(\zeta_n,\sqrt{g_1}):\mathbb Q]=\varphi(n)\epsilon_2(n),{\rm ~where~}
\epsilon_2(n)=\cases{1 &if $D(g_1)|n$;\cr 2 &otherwise.}$$
It then follows (cf. the argument leading up to (\ref{tweededoorsnede})) that
$[\mathbb Q(\zeta_f)\cap \mathbb Q(\zeta_n,\sqrt{g_1}):\mathbb Q(\zeta_{(f,n)})]=\epsilon_2(n)/\epsilon_2([f,n])$.
Thus our field intersection is quadratic over
$\mathbb Q(\zeta_{(f,n)})$ if $D(g_1)\nmid n$ and $D(g_1)|[f,n]$ and equals $\mathbb Q(\zeta_{(f,n)})$ otherwise.
In the latter case $\gamma_0(D(g_1))=1$ and we are done. Next assume that $D(g_1)\nmid n$ and $D(g_1)|[f,n]$.
Write $\gamma=\gamma_0(D(g_1))$. By Lemma \ref{atfundem}
$\gamma$ is a fundamental discriminant and hence by Lemma
\ref{conductor} $\mathbb Q(\sqrt{\gamma})$ is a quadratic
field of conductor $|\gamma|$. Since $\gamma|f$, it follows that $\sqrt{\gamma}\in \mathbb Q(\zeta_f)$.
Likewise, using 
Lemma \ref{atfundem}, it is seen that $\mathbb Q(\sqrt{D(g_1)/\gamma})$ is a field of conductor $|D(g_1)/\gamma|$. Note that the
assumption that $D(g_1)$ divides $[f,n]$ implies that ${D(g_1)\over (f,D(g_1))}|{n\over (f,n)}|n$ and, in case
$\nu_2(f)<\nu_2(D)$, that $\nu_2(n)\ge \nu_2(D(g_1))$. Using this we infer that $\sqrt{D(g_1)/\gamma}\in
\mathbb Q(\zeta_n)$. Since $\sqrt{D(g_1)}\in \mathbb Q(\zeta_n,\sqrt{g_1})$, it follows that
$\sqrt{\gamma} \in \mathbb Q(\zeta_n,\sqrt{g_1})$. The proof will now be completed once we show that
$\sqrt{\gamma}\not\in \mathbb Q(\zeta_{(f,n)})$. Assume that $\sqrt{\gamma}\in \mathbb Q(\zeta_{(f,n)})$. Then
one infers, using that $\gamma$ is a fundamental discriminant and, in case $\nu_2(f)<\nu_2(D)$, that
$\nu_2(n)\ge \nu_2(D)$, that $(f,D(g_1))|(f,n)$. This together with
${D(g_1)\over (f,D(g_1))}|{n\over (f,n)}$ then leads to the conclusion that $D(g_1)|n$, contradicting our
assumption that $D(g_1)\nmid n$. \qed

\begin{Lem}
\label{nieuwedoorsnede}
Let $g_1\ne 0$ be a rational number.
We have
\begin{equation}
\label{reductie}
\mathbb Q(\zeta_n,\zeta_{2n}\sqrt{g_1})=\cases{\mathbb Q(\zeta_n,\sqrt{g_1}) &if $2\nmid n$;\cr
\mathbb Q(\zeta_n,\sqrt{-g_1}) &if $\nu_2(n)=1$;\cr
\mathbb Q(\zeta_n,\sqrt{2g_1}) &if $\nu_2(n)=2$ and $8|D(g_1)$,}
\end{equation}
and in these cases the intersection 
$I=\mathbb Q(\zeta_f)\cap \mathbb Q(\zeta_n,\zeta_{2n}\sqrt{g_1})$ can be determined
using Lemma {\rm \ref{oudedoorsnede}}. Next assume that none of the conditions 
{\rm (\ref{reductie})} are satisfied. Then
$$[I:\mathbb Q(\zeta_{(f,n)})]=\cases{2 &if $D(g_1)|[f,n]$ and $\nu_2(f)>\nu_2(n)$;\cr
1 &otherwise.}$$
Furthermore, 
$$I=\cases{\mathbb Q(\zeta_{(f,n)},\zeta_{2(f,n)}
\sqrt{\gamma(D(g_1)}) &if $D(g_1)|[f,n]$ and $\nu_2(f)>\nu_2(n)$;\cr
\mathbb Q(\zeta_{(f,n)}) &otherwise.}$$
\end{Lem}
{\it Proof}. 
Let us consider the case where $4|n$, not
both $8|D(g_1)$ and $4||n$, and 
such that, moreover, $D(g_1)|[f,n]$ and $\nu_2(f)>\nu_2(n)$. 
As in the proof of Lemma \ref{oudedoorsnede} we infer that $\sqrt{\gamma(D(g_1))}\in \mathbb Q(\zeta_f)$.
Since $\nu_2(f)>\nu_2(n)$, we infer that $\zeta_{2(f,n)}\sqrt{\gamma(D(g_1))}$ is 
in $\mathbb Q(\zeta_f)$.
As in the proof of Lemma \ref{oudedoorsnede} we
infer that $\sqrt{D(g_1)/\gamma(D(g_1))}\in \mathbb Q(\zeta_n)$.
Since $\zeta_{2(f,n)}\sqrt{D(g_1)}\in \mathbb Q(\zeta_n,\zeta_{2n}\sqrt{g_1})$, we infer
that $\zeta_{2(f,n)}\sqrt{\gamma(D(g_1))}\in \mathbb Q(\zeta_n,\zeta_{2n}\sqrt{g_1})$.
It remains to show that $\zeta_{2(f,n)}\sqrt{\gamma(D(g_1))}\not\in \mathbb Q(\zeta_{(f,n)})$.
If the latter element would be contained in $\mathbb Q(\zeta_{(f,n)})$, by Lemma \ref{degree1}
we would have $D(2\gamma(D(g_1)))|(f,n)$ and $8||2(f,n)$. 
Note that these conditions imply that $8|D(2g_1)$ and $8|f$ and
so $8|D(2\gamma(D(g_1)))|(f,n)$, contradicting the assumption that $4||(f,n)$.
Note that
\begin{equation} 
\label{henk} 
[I:\mathbb Q(\zeta_{(f,n)})]=[\mathbb Q(\zeta_f)\cap \mathbb Q(\zeta_n,\zeta_{2n}\sqrt{g_1}):\mathbb Q(\zeta_{(f,n)})]=
{\epsilon_1(n,2n)\over \epsilon_1([f,n],2n)}.
\end{equation}  
The conditions imposed ensure that $\epsilon_1(n,n)=1$. Using (\ref{henk}) and Lemma \ref{degree1} it follows
that
$$[I:\mathbb Q(\zeta_{(f,n)})]=\cases{2 & if $D(g_1)|[f,n]$ and $\nu_2(f)>\nu_2(n)$;\cr
2 & if $D(2g_1)|[f,n]$ and $\nu_2(n)=2$;\cr
1 & otherwise.}$$
Let us assume that $D(2g_1)|[f,n]$ and $\nu_2(n)=2$. Since if $\nu_2(n)=2$, then $8\nmid D(g_1)$ (by assumption), it
follows that $\nu_2(f)>\nu_2(n)$ and so $D(g_1)|[f,n]$. It follows that
$$[I:\mathbb Q(\zeta_{(f,n)})]=\cases{2 & if $D(g_1)|[f,n]$ and $\nu_2(f)<\nu_2(n)$;\cr
1 & otherwise.}$$
This concludes the proof. \qed\\

\noindent Let $K_{n,n}^{\rm ab}$ denote the maximal abelian subfield of $K_{n,n}$. Note that
$\mathbb Q(\zeta_f)\cap K_{n,n}=\mathbb Q(\zeta_f)\cap K_{n,n}^{\rm ab}$. In the
next lemma $K_{n,n}^{\rm ab}$ is determined. Note that for each choice of $g$ precisely one of the seven cases
applies.
\begin{Lem}
\label{gevallen}
Let $g$ and $g_0$ be as in Lemma {\rm \ref{degree}}.
Let $K_{n,n}^{\rm ab}$ denote the maximal abelian subfield of $K_{n,n}$. We have\\
{\rm 1)} $K_{n,n}^{\rm ab}=\mathbb Q(\zeta_n)$ if $g>0$ and $\nu_2(n)\le \nu_2(h)$;\\
{\rm 2)} $K_{n,n}^{\rm ab}=\mathbb Q(\zeta_n,\sqrt{g_0})$ if $g>0$ and $\nu_2(n)\ge \nu_2(h)+1$;\\
{\rm 3)} $K_{n,n}^{\rm ab}=\mathbb Q(\zeta_n,\sqrt{g_0})$ if $g<0$ and $\nu_2(n)\ge \nu_2(h)+2$;\\
{\rm 4)} $K_{n,n}^{\rm ab}=\mathbb Q(\zeta_n,\sqrt{-g_0})$ if $g<0,~\nu_2(n)=1$ and $\nu_2(h)=0$;\\
{\rm 5)} $K_{n,n}^{\rm ab}=\mathbb Q(\zeta_n,\sqrt{2g_0})$ if $g<0,~\nu_2(h)=1,~\nu_2(n)=2$ and $8|D(g_0)$;\\
{\rm 6)} $K_{n,n}^{\rm ab}=\mathbb Q(\zeta_{2n})$ if $g<0$ and $\nu_2(n)\le \nu_2(h)$;\\
{\rm 7)} $K_{n,n}^{\rm ab}=\mathbb Q(\zeta_n,\zeta_{2n}\sqrt{g_0})$ if $g<0$ and $\nu_2(n)=\nu_2(h)+1$ and none
of the previous cases apply.
\end{Lem}
\noindent The latter result can be formulated more compactly, but with regards to computing $\mathbb Q(\zeta_f)\cap K_{n,n}$
the more extended formulation turns out to be more suitable. A more compact reformulation of Lemma \ref{gevallen}
is as follows.
\begin{Prop}
If $g>0$, then
$$K_{n,n}^{\rm ab}=\cases{\mathbb Q(\zeta_n) &if $\nu_2(n)\le \nu_2(h)$;\cr
\mathbb Q(\zeta_n,\sqrt{g_0}) &if $\nu_2(n)>\nu_2(h)$.}$$
If $g<0$, then
$$K_{n,n}^{\rm ab}=\cases{\mathbb Q(\zeta_{2n}) &if $\nu_2(n)\le \nu_2(h)$;\cr
\mathbb Q(\zeta_n,\zeta_{2n}\sqrt{g_0}) &if $\nu_2(n)=\nu_2(h)+1$;\cr
\mathbb Q(\zeta_n,\sqrt{g_0}) &if $\nu_2(n)\ge \nu_2(h)+2$.}$$
\end{Prop}
{\it Proof of Lemma} \ref{gevallen}. For each 
of the seven
cases one easily checks that the explicit field indicated, let us
call it $E$, is abelian and is contained in $K_{n,n}^{\rm ab}$. The idea
is now to apply Lemma \ref{doorisdoor} with $L=E$ and $M=K_{n,n}$.
By (\ref{tweededoorsnede}) it follows that $[\mathbb Q(\zeta_f)\cap K_{n,n}:\mathbb Q]=\varphi((f,n))\epsilon([f,n],n)/\epsilon(n,n)$.
In cases 1 and 6 the degree $[\mathbb Q(\zeta_f)\cap E:\mathbb Q]$ is trivially determined. In the cases
2, 3, 4 and 5 it can be determined by invoking Lemma \ref{oudedoorsnede}. In the remaining case 6 we
apply Lemma \ref{nieuwedoorsnede}. One computes that in each of the seven cases 
$[\mathbb Q(\zeta_f)\cap E:\mathbb Q]=[\mathbb Q(\zeta_f)\cap K_{n,n}:\mathbb Q]$ for every $f$ and
hence the result then follows on invoking Lemma \ref{doorisdoor}.\\
\indent As an example we deal with case 2. One checks that
$\mathbb Q(\zeta_n,\sqrt{g_0})\subseteq K_{n,n}^{\rm ab}$. By Lemma \ref{oudedoorsnede} it follows  that
$$[\mathbb Q(\zeta_f)\cap \mathbb Q(\zeta_n,\sqrt{g_0}):\mathbb Q(\zeta_{(f,n)}]=
\cases{2 & if $D(g_0)\nmid n$ and $D(g_0)|[f,n]$;\cr
1 & otherwise.}$$ On the other hand, by (\ref{tweededoorsnede}) and Lemma \ref{degree},
$$[\mathbb Q(\zeta_f)\cap K_{n,n}^{\rm ab}:\mathbb Q(\zeta_{(f,n)}]={\epsilon([f,n],n)\over \epsilon(n,n)}=
\cases{2 & if $D(g_0)\nmid n$ and $D(g_0)|[f,n]$;\cr
1 & otherwise.}$$
The proof of this case then follows on invoking Lemma \ref{doorisdoor}. \qed\\

\noindent Next we compute the intersection of $\mathbb Q(\zeta_f)$ with the fields given in
Lemma \ref{gevallen}. To this end we require a definition.
\begin{Def}
\label{gammag}
In each of the seven cases as described in 
{\rm Lemma \ref{gevallen}}, we define $\gamma_g(f,n)$ as follows.\\
{\rm 1)} $\gamma_g(f,n)=1$.\\
{\rm 2)} $\gamma_g(f,n)=\gamma_0(D(g_0))$.\\
{\rm 3)} $\gamma_g(f,n)=\gamma_0(D(g_0))$.\\
{\rm 4)} $\gamma_g(f,n)=\gamma_0(D(-g_0))$.\\
{\rm 5)} $\gamma_g(f,n)=\gamma_0(D(2g_0))$.\\
{\rm 6)} $\gamma_g(f,n)=1$.\\
{\rm 7)} $\gamma_g(f,n)=\gamma_1(D(g_0))$.
\end{Def}

\begin{Def}
If $\nu_2(f)>\nu_2(n)\ge 1$ and either we are in case $6$ or we are in case $7$ 
and $D(g_0)|[f,n]$, then we say we are in the exceptional case.
\end{Def}
Note that if $h$ is odd or $g>0$, the exceptional case does not arise.\\
\indent Now the main result of this section can be formulated:
\begin{Thm}
\label{maindegree}
The equality
$\mathbb Q(\zeta_f)\cap K_{n,n}=\mathbb Q(\zeta_{(f,n)},\sqrt{\gamma_g(f,n)})$ holds, unless 
we are in the exceptional case when
we have
$$\mathbb Q(\zeta_f)\cap K_{n,n}=\mathbb Q(\zeta_{(f,n)},\zeta_{2(f,n)}\sqrt{\gamma_g(f,n)}).$$
\end{Thm}
{\it Proof}. Follows on combining Lemma \ref{gevallen} with Lemma \ref{nieuwedoorsnede} and
Lemma \ref{oudedoorsnede}. \qed

\begin{Cor} Suppose that either $g>0$ or $g<0$ and $h$ is odd. Then
$$\mathbb Q(\zeta_f)\cap K_{n,n}
=\cases{\mathbb Q(\zeta_{(f,n)},\sqrt{\gamma_0(D(g))}) &if $\nu_2(n)\ge \nu_2(h)+1$;\cr
\mathbb Q(\zeta_{(f,n)}) &otherwise.}$$
\end{Cor}
{\it Proof}. By Lemma \ref{gevallen} we have 
$$K_{n,n}^{\rm ab}=\cases{\mathbb Q(\zeta_n,\sqrt{g}) & if $\nu_2(n)\ge \nu_2(h)+1$;\cr
\mathbb Q(\zeta_n) & otherwise.}$$
Now invoke Lemma \ref{oudedoorsnede}. \qed

\begin{Lem}
\label{chebbie}
Let $a$ and $f$ be natural numbers that are coprime. Then the density of primes $p$ such
that $p$ splits completely in $K_{n,n}$ and $p\equiv a({\rm mod~}f)$ is zero  
if $a\not \equiv 1({\rm mod~}(f,n))$ and equals 
$${1+\epsilon_3(a,f,n)({\gamma_g(f,n)\over a})\over 2[K_{[f,n],n}:\mathbb Q]}$$
otherwise, where 
$({\cdot \over \cdot})$ 
denotes the Kronecker symbol and 
$$\epsilon_3(a,f,n)=\cases{(-1)^{a-1\over (f,n)} &in the exceptional case;\cr
1 &otherwise.}$$
\end{Lem}
{\it Proof}.
Using Chebotarev's Density Theorem we infer that
the set under consideration has a density and that the density equals $1/[K_{[f,n],n}:\mathbb Q]$ if
the restriction of $\sigma_a$ to $\mathbb Q(\zeta_f)\cap K_{n,n}$ is the identity and zero otherwise.
The latter intersection of fields always contains $\zeta_{(f,n)}$. So if $a\not\equiv 1({\rm mod~}(f,n))$, then
$\sigma_a$ does not leave $\zeta_{(f,n)}$ fixed and the density is zero. So assume that
$a\equiv 1({\rm mod~}(f,n))$. Using Lemma \ref{conductor2} and the
fact that $\gamma_g(f,n)$ is a fundamental discriminant, we infer that $\sigma_a(\sqrt{\gamma_g(f,n)})=({\gamma_g(f,n)\over a})
\sqrt{\gamma_g(f,n)}$. Furthermore, we have $\sigma_a(\zeta_{2(f,n)})=(-1)^{(a-1)/(f,n)}\zeta_{2(f,n)}
=\epsilon_3(a,f,n)\zeta_{2(f,n)}$.
On invoking Theorem \ref{maindegree} the result then follows. \qed\\

\noindent {\tt Remark}. If $\alpha\in \mathbb Q(\zeta_f)\cap K_{n,n}$ and $\alpha\not\in \mathbb Q(\zeta_{(f,n)})$,
then $\sigma_a(\alpha)=\epsilon_3(a,f,n)({\gamma_g(f,n)\over a})\alpha$.
\begin{Cor}
\label{dichtheidomschrijving}
Let $n$ be squarefree. Put $t_d=(t,d^{\infty})$.
The density of primes $p$ such that $p\equiv 1+ta({\rm mod~}dt)$ and $p$ splits completely in $K_{nt,nt}$
equals zero if $(d,n)\nmid a$ or $(1+ta,d)>1$, otherwise it equals
\begin{equation}
\label{transfer}
{1+\epsilon_3(1+ta,dt,nt)({\gamma_g(dt,nt)\over 1+ta})\over 2[K_{[d,n]t,nt}:\mathbb Q]}
={1+\epsilon_3(1+ta,dt_d,nt)({\gamma_g(dt_d,nt)
\over 1+ta})\over 2[K_{[d,n]t,nt}:\mathbb Q]},
\end{equation}
where
$$\epsilon_3(1+ta,dt_d,nt)=\cases{(-1)^{a\over (d,n)} &in the exceptional case;\cr
1 &otherwise.}$$
\end{Cor}
{\it Proof}. This follows from Lemma \ref{chebbie} together with the observation that the two systems of
congruences 
$$\cases{x\equiv 1+ta({\rm mod~}dt)\cr x\equiv 1({\rm mod~}nt)}{\rm ~and~}
\cases{x\equiv 1+ta({\rm mod~}dt_d)\cr x\equiv 1({\rm mod~}nt)}$$
are equivalent. \qed\\

\noindent {\tt Remark}. The equivalence of the system of congruences shows that $[dt_d,nt]=[dt,nt]=[d,n]t$, an
observation that will be needed later (e.g. in the proof of Lemma \ref{DDT}).\\

\noindent The advantage of the formula on the right hand side of (\ref{transfer}) is that the numbers
$\gamma_g(dt_d,nt)$ that occur are restricted to a finite set of divisors, namely those of 
$(4d,D(g_0))$, so this is a very restricted set of numbers. The situation is described more precisely
by the next result (the easy proof of which is left to the reader).
\begin{Lem}
\label{DDT}
Let $D$ be a fundamental discriminant.\\
If $(d,D)=1$, then $\gamma_0(D;dt_d,nt)=1$. Assume that $(d,D)\ne 1$.
Set 
$$\Delta_1(D)=(-1)^{(d_{\rm odd},D)-1\over 2}(d_{\rm odd},D)~{\rm and~}\Delta_2(D)=(-1)^{D/(8d,D)-1\over 2}(8d,D).$$
{\rm 1)} If $2\nmid d$, then 
$$\gamma_0(D;dt_d,nt)=\cases{\Delta_1(D) & if $D\nmid nt$ and $D|[d,n]t$;\cr
1 & otherwise,}$$
where $\Delta_1(D)|d$.\\
{\rm 2)} If $\nu_2(d)\ge \nu_2(D)$, then 
$$\gamma_0(D;dt_d,nt)=\cases{\Delta_2(D) & if $D\nmid nt$ and $D|[d,n]t$;\cr
1 & otherwise,}$$
where $\Delta_2(D)|d$.\\
{\rm 3)} If $1\le \nu_2(d)<\nu_2(D)$, then
$$\gamma_0(D;dt_d,nt)=\cases{\Delta_1(D) & if $\nu_2(t)<\nu_2(D/d)$, $D\nmid nt$ and $D|[d,n]t$;\cr
\Delta_2(D) & if $\nu_2(t)\ge \nu_2(D/d)$, $D\nmid nt$ and $D|[d,n]t$;\cr
1 & otherwise,}$$
where $\Delta_1(D)|d/2$ and $\Delta_2(D)|4d$.
\end{Lem}

\noindent The symbol $({\gamma_g(f,n)\over a})$ appearing in Lemma \ref{chebbie} has a preference
(certainly as $D(g_0)$ becomes large), to be equal to 1. This can be quantified as follows.
\begin{Prop}
\label{symboliseen}
Let $a$ and $f$ be natural numbers that are coprime and assume that $a\equiv 1({\rm mod~}(f,n))$.
Then $({\gamma_g(f,n)\over a})\ne 0$. If $D(g_0)_{\rm odd}\nmid [f,n]$, then
$({\gamma_g(f,n)\over a})=1$. 
\end{Prop}
{\it Proof}. Since by definition $\gamma_g(f,n)|f$ and by assumption $(a,f)=1$, the definition of the Kronecker
symbol implies that $({\gamma_g(f,n)\over a})\ne 0$. If $D(g_0)_{\rm odd}\nmid [f,n]$, then
Definition \ref{gammag} implies that $\gamma_g(f,n)=1$. On noting 
that  $({\gamma_g(f,n)\over a})=({1\over a})=1$, the result follows. \qed\\

\noindent {\tt Remark}. If $h$ is odd the conclusion of Proposition \ref{symboliseen} holds if we make the (weaker)
assumption that
$D(g)\nmid [f,n]$.

\section{A more explicit formula  for $\delta_g(a,d)$}
After all this preparation we can rewrite (\ref{centraal}) in a more explicit form involving 
the Kronecker symbol.
\begin{Thm} 
\label{main}
{\rm (GRH)}. Let $a\ge 1$ and $t_d=(t,d^{\infty})$. Then
$$\delta_g(a,d)=\sum_{t=1\atop (1+ta,d)=1}^{\infty}
\sum_{n=1\atop (n,d)|a}^{\infty}{\mu(n)c_g(1+ta,dt,nt)\over [K_{[d,n]t,nt}:\mathbb Q]},$$
where
$$c_g(1+ta,dt,nt)={1\over 2}\left(1+\epsilon_3(1+ta,dt_d,nt)\left({\gamma_g(dt_d,nt)\over 1+ta}\right)\right),$$
with
$\gamma_g(\cdot,\cdot)$ as in Definition {\rm \ref{gammag}},  $({\cdot\over \cdot})$ the Kronecker symbol 
and where $\epsilon_3(1+ta,dt_d,nt)$
equals $(-1)^{a/(d,n)}$ if $\nu_2(d)>\nu_2(n)$, $\nu_2(nt)\ge 1$ and $K_{nt,nt}^{\rm ab}$ falls
under case $6$ or $7$ of Lemma {\rm \ref{gevallen}} and where  $\epsilon_3(1+ta,dt_d,nt)=1$ otherwise.
\end{Thm}
{\it Proof}. Follows from (\ref{centraal}), the observation that
$c_g(1+ta,dt,nt)/[K_{[d,n]t,nt}:\mathbb Q]$ is the density of primes $p\equiv 1+ta({\rm mod~}dt)$ 
that split completely in $K_{nt,nt}$, Corollary \ref{dichtheidomschrijving} and the definition of
the exceptional case. \qed\\

\noindent {\tt Remark}. If $g>0$, $d$ is odd or $h$ is odd the exceptional case does not arise and 
hence $\epsilon_3(.,.,.)=1$ for the relevant values of $t$ and $n$. \\

\noindent {\tt Example}. Suppose that $h$ is odd and that, furthermore, $d$ is odd or $8|d$. Then, on GRH,
\begin{equation}
\label{laternodig}
\delta_g(a,d)=\delta_g^{(0)}(a,d)-\sum_{t=1\atop {(1+ta,d)=1\atop ({\delta\over 1+ta})=-1}}^{\infty}
\sum_{n=1,~(d,n)|a\atop {D(g)\nmid nt,~D(g)|[d,n]t\atop \nu_2(nt)\ge 1}}^{\infty}{\mu(n)\over [K_{[d,n]t,nt}:\mathbb Q]},
\end{equation}
where 
$$\delta=\cases{\Delta_1(D(g)) &if $d$ is odd;\cr
\Delta_2(D(g)) &if $8|d$,}$$
where $\Delta_1$ and $\Delta_2$ are as in Lemma \ref{DDT} and divide $d$. Alternatively one can write
$$\delta_g(a,d)=\delta_g^{(0)}(a,d)-\sum_{t=1\atop {(1+ta,d)=1\atop ({\delta\over 1+ta})=-1}}^{\infty}
\sum_{n=1,~(n,d)|a\atop \sqrt{\delta}\in K_{nt,nt}}^{\infty}{\mu(n)\over [K_{[d,n]t,nt}:\mathbb Q]}.$$

\noindent {\tt Example}. If $g>0$ and either $d$ is odd or $8|d$, then, on GRH, (\ref{laternodig}) holds with
the condition $\nu_2(nt)\ge 1$ replaced by $\nu_2(nt)\ge \nu_2(h)+1$ and $D(g)$ by $D(g_0)$.\\

\noindent {\tt Example} (occurrence of exceptional case). Suppose that $8|d$, $g<0$ and $\nu_2(h)\ge 2$. Let $t$
and $n$ be such that $(1+ta,d)=1$ and $(n,d)|a$. Let $\delta=\Delta_2(D(g_0))$. Then
$$
2c_g(1+ta,dt,nt)=\cases{1+({\delta\over 1+ta}) & if $\nu_2(nt)\ge \nu_2(h)+2$;\cr
1+(-1)^{a\over (d,n)}({\delta\over 1+ta}) & if $\nu_2(nt)=\nu_2(h)+1$ and $D(g_0)|[d,n]t$;\cr
2 & otherwise.}
$$

\noindent In case $1\le \nu_2(d)\le 2$ more complicated formulae for $\delta_g(a,d)$ arise. For our purposes
it turns out, however, that in order to prove Theorem \ref{favoriet} or part 2 of Theorem \ref{indeks} these
are not needed.

\section{Proof of Theorems \ref{peelof} and \ref{stellingeen}}
{\it Proof of Theorem} \ref{peelof}.
A proof of part 1, based on part 1 of Lemma \ref{graadoprekking} was given in \cite{Moreealleen2}. Similarly part 2 can
be proved using part 2 of Lemma \ref{graadoprekking}. The result for 
$\delta_g^{(0)}$ follows trivially from Lemma \ref{graadoprekking} and for $\delta$
it follows easily from part 1 of Theorem 1 of \cite{Moreealleen}.\qed\\

\noindent {\it Proof of Theorem} \ref{stellingeen}.\\
{\it Proof of the inequality for} $|\delta_g^{(0)}(a,d)-\delta(a,d)|$.  Using Lemma \ref{degree} we infer that
\begin{equation}
\label{vergmetgem}
\delta_g^{(0)}(a,d)
=\sum_{t=1\atop (1+ta,d)=1}^{\infty}\sum_{n=1\atop (n,d)|a}^{\infty}{\mu(n)\over \varphi([d,n]t)nt}+E,
\end{equation}
with 
$$|E|\le \sum_{t=1\atop (1+ta,d)=1}^{\infty}\sum_{n=1\atop (n,d)|a,~D(g)|dnt}^{\infty}{|\mu(n)|\over \varphi([d,n]t)nt}.$$
Let us denote the latter double sum by $E_1$.
Note that
\begin{eqnarray}
\label{hulpbijeersteschatting}
E_1 &\le & \sum_{t=1}^{\infty}\sum_{n=1\atop D(g)|dnt}^{\infty}{|\mu(n)|\over \varphi(nt)nt}
=\sum_{D_1|v}{2^{\omega(v)}\over \varphi(v) v}\cr
&=&\sum_{v_1=1}^{\infty}{2^{\omega(D_1v_1)}\over \varphi(D_1v_1)D_1v_1}\le 
{2^{\omega(D_1)}\over \varphi(D_1)D_1}\sum_{v=1}^{\infty}
{2^{\omega(v)}\over \varphi(v) v}<{2^{\omega(D_1)+2}\over \varphi(D_1)D_1}.
\end{eqnarray}
By Theorem 3 of \cite{Moreeaverage} the second double sum in (\ref{vergmetgem}) equals $\delta(a,d)$. \qed\\

\noindent {\tt Remark}. In case $(a,d)=1$, the term $\varphi(D_1)$ in the first inequality of Theorem 
\ref{stellingeen} can be replaced by $\varphi(dD_1)$.\\ 

\noindent {\it Proof of the inequality for} $|\delta_g(a,d)-\delta(a,d)|$. Using Theorem \ref{main} and the 
remark following Proposition \ref{symboliseen} we have
$\delta_g(a,d)=\delta_g^{(0)}(a,d)+E_2$,
with
$$|E_2|\le \sum_{t=1\atop (1+ta,d)=1}^{\infty}
\sum_{n=1\atop (n,d)|a,~D(g)|[d,n]t}^{\infty}{|\mu(n)|\over [K_{[d,n]t,nt}:\mathbb Q]}.$$
Using Lemma \ref{degree} we infer that
$$|E_2|\le 2\sum_{t=1\atop (1+ta,d)=1}^{\infty}\sum_{n=1\atop (n,d)|a,~D(g)|dnt}^{\infty}{|\mu(n)|\over \varphi([d,n]t)nt}
=2E_1.$$
Using (\ref{hulpbijeersteschatting}), we then infer
that $|E_2|<2^{\omega(D_1)+3}/(\varphi(D_1)D_1)$. From the latter estimate, (\ref{hulpbijeersteschatting}) 
and the inequality for $\delta_g^{(0)}(a,d)$, the result then follows. \qed\\ 

\noindent {\it Proof of inequality for} $|\rho_g(a,d)-\rho(a,d)|$. For this, see 
the proof of Proposition 5 of \cite{Moreeaverage}. \qed

\section{The index reconsidered}
In this section we compare $\rho_g(a,d)$ with its average value $\rho(a,d)$. We will use the following
two lemmas.
\begin{Lem}
\label{wirsinggevolg} {\rm \cite{Moreealleen2}}.
Let $d\ge 3$. The number of integers $1\le g\le x$ such that $D(g)$ has
no prime divisor $p$ with $p\equiv 1({\rm mod~}d)$ is
$O_d(x\log^{-1/\varphi(d)}x)$. The same assertion holds with
$D(g)$ replaced by $D(-g)$.
\end{Lem}
The following simple result is also needed.
\begin{Lem}
\label{hijisnul}
If $p|v$ with $p\equiv 1({\rm mod~}d)$ and $p\nmid m$, then
$$\sum_{t\equiv a({\rm mod~}d)\atop t|v,~({v\over t},m)|r}\mu({v\over t})=0.$$
\end{Lem}
{\it Proof}. Let $p^e||v$. Put $v_1=v/p^e$. Then
$$\sum_{t\equiv a({\rm mod~}d)\atop t|v,~({v\over t},m)=\alpha}\mu({v\over t})=
\sum_{j=0}^e \sum_{t\equiv a({\rm mod~}d)\atop t|v_1,~({v_1\over t},m)=\alpha}\mu(p^j{v_1\over t})
=\sum_{j=0}^1 (-1)^j \sum_{t\equiv a({\rm mod~}d)\atop t|v_1,~({v_1\over t},m)=\alpha}
\mu({v_1\over t})=0.$$
\begin{Prop} 
\label{bijnagelijkrho}
{\rm (GRH)}. Let $\alpha=(a,d)$ and $d_1=d/\alpha$. Suppose that $g\in G$ and $D(g)$ has
a prime divisor $p$ with $p\equiv 1({\rm mod~}d_1)$ and $p\nmid \alpha$, then $\rho_g(a,d)=\rho(a,d)$.
\end{Prop}
{\it Proof}. On GRH we have (see \cite{Moreealleen}) 
$$\rho_g(a,d)=\sum_{t\equiv a({\rm mod~}d)}\sum_{v=1}^{\infty}
{\mu(v)\over [K_{vt,vt}:\mathbb Q]}
=\sum_{w=1}^{\infty}{\sum_{t|w,~t\equiv a({\rm mod~}d)}\mu({w\over t})
\over [K_{w,w}:\mathbb Q]}.$$
By \cite{Moreeaverage} we have
$$\rho(a,d)=\sum_{w=1}^{\infty}{\sum_{t|w,~t\equiv a({\rm mod~}d)}\mu({w\over t})
\over w\varphi(w)}.$$
Let $n_1$ be as in Lemma  
\ref{degree} (note that $n_1=[2,D(g)]$).
The assumption that $g\in G$ ensures,
by Lemma \ref{degree}, that
\begin{eqnarray}
\rho_g(a,d)&=&\rho(a,d)+\sum_{n_1|w}{\sum_{t|w,~t\equiv a({\rm mod~}d)}\mu({w\over t})
\over w\varphi(w)}.\nonumber\\
&=&\rho(a,d)+\sum_{{n_1\over (n_1,\alpha)}|w}{\sum_{t|w,~t\equiv a/\alpha({\rm mod~}d_1)}\mu({w\over t})
\over w\alpha\varphi(w\alpha)}.\nonumber
\end{eqnarray}
If $p|D(g)$, $p\equiv 1({\rm mod~}d_1)$ and $p\nmid \alpha$, then $p|{n_1\over (n_1,\alpha)}$ and so $p$ divides every
$w$ occurring in the latter sum and hence, by Lemma \ref{hijisnul}, this
sum will equal zero. \qed

\begin{Prop} 
\label{rhogisrho}
{\rm (GRH)}. Let $\alpha=(a,d)$ and
$d_1=d/\alpha$. For almost all integers $g$ we have $\rho_g(a,d)=\rho(a,d)$.
More precisely, if $d_1\le 2$, then there are at most $O_d(\sqrt{x})$ integers $g$ with
$|g|\le x$ such that $\rho_g(a,d)\ne \rho_g(a,d)$. If $d_1\ge 3$, then
there are
at most $O_d(x {\log^{-1/\varphi(d_1)}x})$ integers $g$ with $|g|\le x$ such that, moreover,
$\rho_g(a,d)\ne \rho(a,d)$.
\end{Prop}
{\it Proof}. In case $d_1\le 2$, then $\rho_g(a,d)=\rho(a,d)$ for those $g$ for
which  $g\in G$ and $D(g)$ is not in a certain finite set by the previous proposition. The set of remaining $g$
is of size $O_d(\sqrt{x})$. (Note that there are at most $O(\sqrt{x})$ integers $|g|\le x$ that are
not in $G$ and that, for a fixed 
integer $g_1$, there at most $O(\sqrt{x})$ integers $g_2$ satisfying $D(g_2)=D(g_1)$.) If $d_1\ge 3$ we invoke
Lemma \ref{wirsinggevolg} with $d=d_1$. \qed\\

\noindent In case $d|2(a,d)$, that is $d_1\le 2$, the latter result can be improved.
Then $\rho(a,d)$ is a rational number and, more precisely, we have:
\begin{Lem}
\label{null}
{\rm  \cite{Moreeaverage}}. Let $d\ge 1$. Then
$$\rho(0,d)={1\over d\varphi(d)} {\rm ~and~}
\rho(d,2d)=\cases{\rho(0,2d) &if $d$ is odd;\cr
3\rho(0,2d) &if $d$ is even.}$$
\end{Lem}
In this case we can both improve the error term and drop the
GRH assumption in Proposition \ref{rhogisrho}.
\begin{Prop}
\label{scherper}
We have at most $O_d(\sqrt{x})$ integers $g$ with $|g|\le x$ such that
$\rho_g(0,d)\ne \rho(0,d)$. Similarly we have at most 
$O_d(\sqrt{x})$ integers $g$ with $|g|\le x$ such that
$\rho_g(d,2d)\ne \rho(d,2d)$.
\end{Prop}
{\it Proof}. It is not difficult to show that
$\rho_g(0,d)=1/[K_{d,d}:\mathbb Q]$ and hence
$\rho_g(d,2d)=1/[K_{d,d}:\mathbb Q]-
1/[K_{2d,2d}:\mathbb Q]$. 
Using Lemma \ref{degree}, we infer
that
there are at most finitely many squarefree integers $g$ for which
$\rho_g(0,d)\ne 1/d\varphi(d)=\rho(0,d)$ (where the equality
is a consequence of Proposition \ref{null}) and that, similarly, there
are at most finitely many squarefree integers $g$ for which
$\rho_g(0,d)\ne 1/d\varphi(d)-1/2d\varphi(2d)=\rho(0,2d)$ (where
the equality is a consequence of Proposition \ref{null}). From
this the result easily follows. \qed\\

\noindent Theorem \ref{peelof} only seems to have an analog for $\rho_g(a,d)$ in case $d|a$:
\begin{Prop}
\label{nietanaloog}
If $q$ is an odd prime dividing $d_1$, then $\rho_g(0,qd_1)=\rho_g(0,d_1)/q^2$ (on GRH)
and $\rho(0,qd_1)=\rho(0,d_1)/q^2$.
\end{Prop}
{\it Proof}. The first equality is an immediate consequence of (\ref{pappastelling}) and
$[K_{vqd_1,vqd_1}:\mathbb Q]=q^2[K_{vd_1,vd_1}:\mathbb Q]$ which holds by Lemma \ref{degree}
for arbitrary $v$. The second equality is a consequence of Lemma \ref{null}. \qed

\section{Proof of Theorem \ref{favoriet}}
The idea of the proof of Theorem \ref{favoriet} is to show that if $g\in G$ and $D(g)$ contains a prime divisor
$p$ with $p\equiv 1({\rm mod~}k_1(d))$, $p\mid d$, then $*_g(a,d)=*(a,d)$. Lemma \ref{wirsinggevolg} then allows us to
finish the proof. Thus, in outline, the proof is similar to that of Proposition \ref{rhogisrho} in the
previous section. Note that in case $*=\rho$ we are done by Proposition \ref{rhogisrho}.
\begin{Lem}
\label{bijnagelijk}
Let $d\ge 2$. Suppose that $g\in G$ and that $D(g)$ has a prime divisor $p$ with
$p\equiv 1({\rm mod~}d/(a,d))$ and $p\nmid d$, then $\delta_g^{(0)}(a,d)=\delta(a,d)$.
\end{Lem}
{\it Proof}. 
By the definition of $\delta_g^{(0)}(a,d)$ and Lemma \ref{degree}, we find that
$$\delta_g^{(0)}(a,d)=\sum_{t=1\atop (1+ta,d)=1}^{\infty}\sum_{n=1\atop (n,d)|a}^{\infty}{\mu(n)\over \varphi([d,n]t)nt} 
+\sum_{t=1\atop (1+ta,d)=1}^{\infty}\sum_{n=1,~(n,d)|a\atop n_{d/(d,n)}|[d,n]t}^{\infty}{\mu(n)\over \varphi([d,n]t)nt}.$$
In the first double sum we recognize $\delta(a,d)$. Let us denote the second double sum by $E_3$. By the remark
following Lemma \ref{degree} we have $n_{d/(d,n)}=[2^{\nu_2(d/(d,n))+1},D(g)]$. Note that $n_{d/(d,n)}|[d,n]t$
iff $2|nt$ and $D(g)|[d,n]t$.
Following the argument of Proposition 7 of \cite{Moreeaverage} and dragging along the conditions $2|nt$ and $D(g)|[d,n]t$
we obtain
$$E_3=\sum_{\alpha|(a,d)}{\mu(\alpha)\over \alpha}\sum_{t=1\atop (1+ta,d)=1}^{\infty}
\sum_{m=1,~(m,d)=1\atop D(g)|dmt,~2|\alpha mt}^{\infty}{\mu(m)\over \varphi(dmt)mt}.$$
Making the substitution $mt=v$ we obtain
$$E_3=\sum_{\alpha|(a,d)}{\mu(\alpha)\over \alpha}\sum_{v=1\atop D(g)|dv,~2|\alpha v}^{\infty}{\sum_{t|v,~(1+ta,d)=1,~({v\over t},d)=1}\mu({v\over t})
\over \varphi(dv)v}.$$
If $p|D(g)$, $p\equiv 1({\rm mod~}d/(a,d))$ and $p\nmid d$, then $p$ divides every $v$ occurring in the latter sum
and hence, by Lemma \ref{hijisnul}, we deduce that $E_3=0$ (note that the set of integers $t$ satisfying
$(1+ta,d)=1$ is an union of arithmetic progressions of modulus $d/(a,d)$). \qed\\

\noindent {\tt Remark}. This lemma when combined with Theorem \ref{peelof} shows, on GRH, that if $g\in G$, $D(g)$ has a 
prime divisor $p$ with $p\equiv 1({\rm mod~}k_2(d))$ and $p\nmid d$, then $\delta_g^{(0)}(a,d)=\delta(a,d)$.\\

\noindent The next few results are concerned with when the equality $\delta_g(a,d)=\delta_g^{(0)}(a,d)$ holds.
\begin{Prop} {\rm (GRH)}. 
\label{other}
Suppose that $d$ is odd and $(d,D(g_0))=1$. Then $\delta_g(a,d)=\delta_g^{(0)}(a,d)$.
\end{Prop}
{\it Proof}. The condtions  $(d,D(g_0))=1$ and $d$ is odd ensure that $\gamma_g(dt_d,nt)=1$, respectively
$\epsilon_3(1+ta,dt_d,nt)=1$ for all $n$ and $t$ (cf. the remark 
following Theorem \ref{main}). Now apply Theorem \ref{main}. \qed

\begin{Lem}
\label{bijnagelijktwee} 
{\rm (GRH)}
Suppose that $g\in G$ and that $D(g)$ has a prime divisor $p$ with
$p\equiv 1({\rm mod~}k_1(d))$ and $p\nmid d$, then $\delta_g(a,d)=\delta_g^{(0)}(a,d)$.
\end{Lem}
{\it Proof}. Let us first consider the case where $d$ is odd or $8|d$.
It follows from Theorem \ref{peelof} that if $\delta_g(a,k_1(d))=\delta_g^{(0)}(a,k_1(d))$, then $\delta_g(a,d)=
\delta_g^{(0)}(a,d)$ and thus, w.l.o.g., we may assume that $k_1(d)=d$.
By
(\ref{laternodig}) we have
$$\delta_g(a,d)=\delta_g^{(0)}(a,d)-\sum_{(1+ta,d)=1\atop ({\delta\over 1+ta})=-1}
\sum_{D(g)\nmid nt,~D(g)|[d,n]t\atop {\nu_2(nt)\ge \nu_2(h)+1\atop (d,n)|a}}{\mu(n)\over [K_{[d,n]t,nt}:\mathbb Q]},$$
for some $\delta|d$.
Since $\delta$ is a fundamental discriminant dividing $d$ it follows from the reciprocity law for the
Kronecker symbol that the sum over $t$ is a sum over a certain set of arithmetic progressions modulo $d$.
Let us consider the inner sum where $t$ runs over one such arithmetic progression, say $t\ge 1$ and
$t\equiv b({\rm mod~}d)$. Setting $(d,n)=\alpha$ and $nt=v$ we obtain for the inner sum:
$$\sum_{\alpha|(a,d)}\sum_{D(g)\nmid v, D(g)|{d\over \alpha}v\atop \nu_2(v)\ge \nu_2(h)+1}{\sum_{t\equiv b({\rm mod~}d),~({v\over t},d)=\alpha}\mu(v/t)\over
[K_{{d\over \alpha}v,v}:\mathbb Q]}.$$
Now the assumption on $D(g)$ implies, by Lemma \ref{hijisnul}, that for all $v$ occurring in the sum the inner
sum is zero. Hence we infer that $\delta_g(a,d)=\delta_g^{(0)}(a,d)$.\\
\indent In case $d\equiv 4({\rm mod~}8)$ we note that $\delta_g(a,d)=\delta_g(a,2d)+\delta_g(a+d,2d)$ and
 $\delta_g^{(0)}(a,d)=\delta_g^{(0)}(a,2d)+\delta_g^{(0)}(a+d,2d)$. Hence it suffices to show that
 $\delta_g(a,2d)=\delta_g^{(0)}(a,2d)$ and $\delta_g(a+d,2d)=\delta_g^{(0)}(a+d,2d)$. Since $8|2d$ and
 $k_1(2d)=k_1(d)$ we are reduced to the previous case. The remaining case where $2||d$ can be dealt with
 similarly. \qed

\begin{Thm} 
\label{dgisd}
{\rm (GRH)}.  For almost all integers $g$ we have $*_g(a,d)=*(a,d)$, with $*$ is $\delta$ or $\delta^{(0)}$.
More precisely, there are
at most 
$O_d(x\log^{-1/\varphi(k_2(d))}x)$ 
integers $g$ with $|g|\le x$ such that 
$\delta_g^{(0)}(a,d)\ne \delta(a,d)$. There are at most 
$O_d(x \log^{-1/\varphi(k_1(d))}x)$ integers $g$ with $|g|\le x$ such that 
$\delta_g(a,d)\ne \delta(a,d)$.
\end{Thm}
{\it Proof}. This follows from the remark following Lemma \ref{bijnagelijk}, Lemma \ref{bijnagelijktwee} and
Lemma \ref{wirsinggevolg} (cf. the proof of Proposition \ref{rhogisrho}). \qed\\
 
\noindent {\tt Remark}. For various choices of $a$ and $d$, the terms $\varphi(k_2(d))$ and $\varphi(k_1(d))$ 
can be sometimes replaced by smaller numbers, resulting in a smaller error term.\\

\noindent {\it Proof of Theorem} \ref{favoriet}. Combine Theorem \ref{dgisd} and Proposition 
\ref{rhogisrho}. \qed

\section{The proof of Theorem \ref{indeks}}
\noindent For any Dirichlet character $\chi$, we let $h_{\chi}$ denote the
Dirichlet convolution of $\chi$ and $\mu$. For properties of $h_{\chi}$
the reader is referred to \cite{Moreealleen}. 
\begin{Lem}
\label{paulmineen}
Let $\alpha,\beta,d,d_1$ and $r$ be arbitrary positive integers with $(\alpha,\beta)=1$.
Write $v=v_1v_d$ with $v_1=v/(v,d^{\infty})$ and $v_d=(v,d^{\infty})$.
Then
$$\sum_{t=1\atop t\equiv \alpha({\rm mod~}\beta)}^{\infty}
\sum_{n=1\atop (n,d)=1}^{\infty}{\mu(n)\over [K_{d_1rnt,rnt}:\mathbb Q]}
={1\over \varphi(\beta)}\sum_{\chi\in G_{\beta}}{\overline{\chi(\alpha)}}\sum_{v=1}^{\infty}
{\chi(v_d)h_{\chi}(v_1)\over [K_{d_1rv,rv}:\mathbb Q]}.$$
\end{Lem}
{\it Proof}. We have
\begin{eqnarray}
\varphi(\beta)\sum_{t=1\atop t\equiv \alpha({\rm mod~}\beta)}^{\infty}
\sum_{n=1\atop (n,d)=1}^{\infty}{\mu(n)\over [K_{d_1rnt,rnt}:\mathbb Q]}
&=&\sum_{v=1}^{\infty}{\sum_{t\equiv \alpha({\rm mod~}\beta),~t|v,~(v/t,d)=1}\mu(v/t)\over
[K_{d_1vr,vr}:\mathbb Q]}\nonumber\\
&=&\sum_{\chi\in G_{\beta}}{\overline{\chi(\alpha)}}\sum_{v=1}^{\infty}
{\sum_{t|v,~(v/t,d)=1}\chi(t)\mu(v/t)\over [K_{d_1rv,rv}:\mathbb Q]}\nonumber\\
&=&\sum_{\chi\in G_{\beta}}{\overline{\chi(\alpha)}}\sum_{\delta|d^{\infty}}
\chi(\delta)\sum_{v=1\atop (v,d)=1}^{\infty}{h_{\chi}(v)\over 
[K_{d_1\delta rv,\delta rv}:\mathbb Q]}\nonumber\\
&=&\sum_{\chi\in G_{\beta}}{\overline{\chi(\alpha)}}\sum_{v=1}^{\infty}
{\chi(v_d)h_{\chi}(v_1)\over [K_{d_1rv,rv}:\mathbb Q]}.\nonumber
\end{eqnarray}
This concludes the proof. \qed\\

\noindent Note that the inner three sums in the next result are finite sums. The complicated
nature of the result suggests that the arithmetic function $h_{\chi}(v)$ might not be
the proper one to work with. However, in the analysis of $d=3,4$ \cite{Moreealleen} and $d$ is
prime \cite{Moreealleen2}, this function turned out to be of crucial importance.
\begin{Thm} 
\label{nulexpl}
Let $a$ and $d\ge 1$ be arbitrary integers. Then
$$\delta_g^{(0)}(a,d)=
\sum_{0<t_1\le d\atop (1+t_1a,d)=1}
\sum_{\chi\in G_{d\over (t_1,d)}}{\overline{\chi({t_1\over (t_1,d)})}\over 
\varphi({d\over (t_1,d)})}
\sum_{\alpha|(a,d)}\mu(\alpha)
\sum_{v=1}^{\infty}
{\chi(v_d)h_{\chi}(v_1)\over [K_{d(t_1,d)v,d(t_1,d)\alpha v}:\mathbb Q]}.$$
\end{Thm}
{\it Proof}. Note that
\begin{eqnarray}
\sum_{n=1\atop (n,d)|a}^{\infty}{\mu(n)\over [K_{[d,n]t,nt}:\mathbb Q]}
&=&\sum_{\alpha|(a,d)}\sum_{n=1\atop (n,d)=\alpha}^{\infty}{\mu(n)\over
[K_{dnt/(d,n),nt}:\mathbb Q]}\nonumber\\
&=&\sum_{\alpha|(a,d)}\mu(\alpha)
\sum_{n=1\atop (n,d)=1}^{\infty}{\mu(n)\over [K_{dnt,\alpha nt}:\mathbb Q]},\nonumber
\end{eqnarray}
cf. the proof of Part 2 of Theorem 1 \cite{Moreeaverage}.
On invoking Lemma \ref{paulmineen} the proof is then easily completed. \qed

\begin{Lem} 
\label{nulletje}
We have $\delta_g^{(0)}(a,d)=\sum_{\chi\in G_{d}}c_{\chi}A_{\chi}$, where
$c_{\chi} \in \mathbb Q(\zeta_{o_{\chi}})$ can be explicitly computed.
\end{Lem}
The proof of this result can be easily inferred from the next lemma in combination
with Theorem \ref{nulexpl} and the observation that if a Dirichlet character
$\chi$ mod $d$ is the lift of a Dirichlet character $\chi'$ mod $d_1$, then
$A_{\chi}/A_{\chi'}$ is rational. 
\begin{Lem}
\label{verjaardag}
Let $v_1$ and $v_d$ be as in Lemma {\rm \ref{paulmineen}}.
We have $$\sum_{v=1}^{\infty}
{\chi(v_d)h_{\chi}(v_1)\over [K_{d_1rv,rv}:\mathbb Q]}=c_{\chi}A_{\chi},$$
where $c_{\chi}\in \mathbb Q(\zeta_{o_{\chi}})$ can be explicitly computed.
\end{Lem}
{\it Proof}. Using Lemmas 10 and 11 of \cite{Moreealleen} it is easily inferred
that $$\sum_{(r,v)=1}{h_{\chi}(v)\over [K_{sv,v}:\mathbb Q]}=c_{\chi}'A_{\chi},$$ where
$c_{\chi}'$ can be effectively computed. The proof of the present lemma is a variation
of this. \qed\\

\noindent {\it Proof of Theorem} \ref{indeks}.\\
1) Follows from Lemma \ref{nulletje} and Theorem \ref{peelof}.\\
2) (Sketch). The same trick as in the proof of Lemma \ref{bijnagelijktwee} can be applied to reduce to the case
where $d$ is odd or $8|d$. In these cases comparatively easy explicit formulae for $\delta_g(a,d)$
can be written down. Each of them can be expressed in terms of a function very closely related to
$h_{\chi}(v)$, this being possible due to the fact that the $t$ arising in these formulae run over
an union of arithmetic progressions modulo a divisor of $4d$. We then get a finite combination
of sums not unlike those in Lemma \ref{verjaardag} that, like the sums in Lemma \ref{verjaardag}
can be expressed in terms of $A_{\chi}$. As the amount of work involved is quite considerable and
no new ideas are involved, I suppress the details here.\\
3) This is Theorem 5 of \cite{Moreealleen}. \qed

\section{Equidistribution results}

\noindent Let $N_g(a_1,d_1;a_2,d_2)(x)$ denote the number of primes $p\le x$ with $p\equiv a_1({\rm mod~}d_1)$ such
that $\nu_p(g)=0$ and the order of $g$ modulo $p$ is congruent to $a_2({\rm mod~}d_2)$. Although $N_g(a,d)(x)$ does
not seem to have equidistribution properties, results in this direction can be proved for $N_g(a_1,d_1;a_2,d_2)(x)$. 
In \cite{Moreealleen2} it is established that, on GRH, $N_g(a_1,d_1;a_2,d_2)$ has a density $\delta_g(a_1,d_1;a_2,d_2)$.\\
\indent Recall that $k(d)$ denotes the squarefree kernel of $d$. The result below generalizes results 
in \cite{Moreealleen} and \cite{Moreealleen2}.
\begin{Thm} {\rm (GRH)}. Suppose that $(a,d)=(b,d)=1$.\\
{\rm 1)} If $d$ is odd, then
$\rho_g(1,k(d);a,d)=\rho_g(1,k(d);b,d)$.\\
{\rm 2)} If $d$ is even, then
$\rho_g(1,2k(d);a,d)=\rho_g(1,2k(d);b,d)$.
\end{Thm}
{\it Proof}. A variation of Theorem \ref{main} yields that
\begin{equation}
\label{transfervergelijking}
\rho_g(1,k(d);a,d)=\sum_{k(d)|t}\sum_{n=1\atop (n,d)=1}^{\infty}
{\mu(n)\epsilon_4(a,n,t)\over  [K_{[d,n]t,nt}:\mathbb Q]},
\end{equation}
where $\epsilon_4(a,n,t)=c_g(1+ta,dt,nt)=(1+\epsilon_3(1+ta,dt_d,nt)({\gamma_g(dt_d,nt)
\over 1+ta}))/2\in \{0,1\}$.\\
1) The assumption that $d$ be odd, implies that $dt_d$ is odd and hence the exceptional case never arises
and we have $\epsilon_3(1+ta,dt_d,nt)=1$ for the values of $n$ and $t$ we sum over. Using that all quadratic 
subfields of $\mathbb Q(\zeta_{dt_d})$ are contained in $\mathbb Q(\zeta_{k(d)})$ and since $k(d)|t$, it
follows that $\sigma_{1+ta}(\zeta_{k(d)})=\zeta_{k(d)}$ and so certainly
$\sigma_{1+ta}(\sqrt{\gamma_g(dt_d,nt)})=\sqrt{\gamma_g(dt_d,nt)}$ and hence
$({\gamma_g(dt_d,nt)/1+ta})=1$ (alternatively one can use quadratic reciprocity to infer this).
It follows that $\epsilon_4(a,n,t)=1$ for all the values of $n$ and $t$ we sum over.
Similarly we have $\epsilon_4(b,n,t)=1$
and it thus follows that
$$\rho_g(1,k(d);a,d)=\rho_g(1,k(d);b,d)=\sum_{k(d)|t}\sum_{n=1\atop (n,d)=1}^{\infty}
{\mu(n)\over  [K_{[d,n]t,nt}:\mathbb Q]}.$$
2) In this case the analog of (\ref{transfervergelijking}) holds, but with the condition $k(d)|t$ 
replaced by $2k(d)|t$. Now the exceptional case may arise. 
In that case
we have $\epsilon_3(1+ta,dt_d,nt)=(-1)^{ta/t_d}=-1$ and, similarly, $\epsilon_3(1+tb,dt_d,nt)=-1$.
Using that all quadratic 
subfields of $\mathbb Q(\zeta_{dt_d})$ are contained in $\mathbb Q(\zeta_{4k(d)})$ and since $4k(d)|t(a-b)$ 
(since $2k(d)|t$ and $a-b$ is even), it
follows that $\sigma_{1+ta}(\zeta_{4k(d)})=\sigma_{1+tb}(\zeta_{4k(d)})$ and so certainly
$\sigma_{1+ta}(\sqrt{\gamma_g(dt_d,nt)})$ is equal to $\sigma_{1+tb}(\sqrt{\gamma_g(dt_d,nt)})$ and hence
$({\gamma_g(dt_d,nt)/1+ta})=({\gamma_g(dt_d,nt)/1+tb})$.
It follows that $\epsilon_4(a,n,t)=\epsilon_4(b,n,t)$ for all the values of $n$ and $t$ we sum over
in (\ref{transfervergelijking}) and hence
$\rho_g(1,2k(d);a,d)=\rho_g(1,2k(d);b,d)$. \qed


\begin{thebibliography}{99999}
\bibitem[CM]{CM} K. Chinen and L. Murata, On a distribution property of the residual 
order of $a{\rm mod~}p)$, I, II {\it J. Number Theory} {\bf 105} (2004), 60--81, 82-100.
\bibitem[GP]{GP} A.J. Granville and P.A.B. Pleasants, Aurefeuillian factorization
revisited, preprint. 
\bibitem[H]{Hooley} {C. Hooley}, {Artin's conjecture
for primitive roots}, {\it J. Reine Angew. Math.} {\bf 225} (1967),
209--220.
\bibitem[L]{Lenstra} {H. W. Lenstra, jr.}, {On Artin's conjecture and Euclid's
algorithm in global fields}, {\it Invent. Math.} {\bf 42} (1977), 201--224.
\bibitem[M-Fi]{PFibo} P. Moree, Convoluted convolved Fibonacci numbers, 
{\it J. Integer Seq.} {\bf 7} (2004), Article 04.2.2, 16 pp.  (electronic) 
\bibitem[M-Av]{Moreeaverage} {P. Moree}, On the average
number of elements in a finite field with
order or index in a prescribed residue class, arXiv:math.NT/0212220, {\it Finite
Fields Appl.} to appear.
\bibitem[M-I]{Moreealleen} {P. Moree}, On the distribution of the order and index of
$g({\rm mod~}p)$ over residue classes I, math.NT/0211259, submitted.
\bibitem[M-II]{Moreealleen2} {P. Moree}, On the distribution of the order and index of
$g({\rm mod~}p)$ over residue classes II, math.NT/0404339, submitted.
\bibitem[O]{Odoni} R.W.K. Odoni,
A conjecture of Krishnamurthy on decimal periods and some allied problems,
{\it J. Number Theory} {\bf 13} (1981), 303--319.
\bibitem[P]{Pappalardi} F. Pappalardi, On Hooley's theorem with weights, 
Number theory, II (Rome, 1995), {\it Rend. Sem. Mat. Univ. Politec. Torino}
{\bf 53} (1995), 375--388. 
\bibitem[Wa]{Wagstaff} {S.S. Wagstaff, jr.},
Pseudoprimes and a generalization of Artin's conjecture, {\it Acta Arith.}
{\bf 41} (1982), 141--150.
\bibitem[Wi1]{Wiertelak1} K. Wiertelak, On the density of some sets of primes,
IV, {\it Acta Arith.} {\bf 43} (1984), 177--190.
\bibitem[Wi2]{Wiertelak2} K. Wiertelak,
On the density of some sets of primes
$p$, for which $n|{\rm ord}_p(a)$,
{\it Funct. Approx. Comment. Math.} {\bf 28}
(2000), 237--241.
\end{thebibliography}
\end{document}